\documentstyle[amssymb,amsmath,amsthm,amscd,graphicx,12pt]{article}

\title{Functoriality of Isovariant Homotopy Classification}
\author{Sylvain Cappell\thanks{Research was partially supported by an NSF grant} \\
Courant Institute, New York University \\
\\
Shmuel Weinberger\thanks{Research was partially supported by an NSF grant} \\
University of Chicago \\
\\
Min Yan\thanks{Research was supported by Hong Kong Research Grant Council General Reseaqrch Fund HKUST605005} \\ 
Hong Kong University of Science and Technology}
\date{}

\newcommand{\s}{\hspace{2mm}}

\newcommand{\ra}{\rightarrow}

\newcommand{\sub}{\subset}
\newcommand{\lra}{\longrightarrow}

\newcommand{\pa}{\partial}

\newcommand{\mc}{\mathcal}
\newcommand{\bb}{\mathbb}
\newcommand{\rel}{\text{rel }}

\newtheorem{thm}{Theorem}

\newtheorem{prop}[thm]{Proposition}

\theoremstyle{definition}

\numberwithin{equation}{section}

\begin{document}

\maketitle


\begin{abstract}
It is a deep fact that the homotopy classification of topological manifolds is convariantly functorial. In other words, a map $M\to N$ of topological manifolds naturally induces a map $S(M)\to S(N)$ of their structure sets. We extend the fact to the isovariant structure set ${\mc S}_G(M,\text{\rm rel } M_s)$ of $G$-equivariant topological manifolds isovariantly homotopy equivalent to $M$ and restricts to homormorphism on the singular part $M_s$, consisting of those points fixed by some non-trivial elements of $G$. We further explain that the structure set ${\mc S}_G(M,\text{\rm rel } M_s)$ is the fibre of the assembly map for the generalized homology theory with the $L$-spectrum as the coefficient. This relates our result to the Farrell-Jones Conjecture for $L$-theory. 
\end{abstract}

\section{Introduction}

The natural problems regarding group actions are usually {\em equivariant}. The maps involved in such problems are only required to preserve the group actions. For example, any map into the single point space is equivariant. This poses some technical difficulties for studying such problems directly. The classical approach often assumes that the fixed points of smaller subgroups have dimensions more than twice of the fixed points of bigger subgroups (big gap hypothesis). The assumption allows enough room to carry out classical manipulations on fixed points of various isotropy subgroups without mutual interference. An alternative method is to consider {\em isovariant} problems, in which maps not only preserves the group actions but also the isotropy groups. Thus a map into the single point space is isovariant only when the action on the source space is trivial. The advantage of the approach is that the noninterference between fixed points of various isotropy subgroups is built into the definition, so that the big gap hypothesis is no longer needed. The obvious drawback is the artificial nature of the isovariant condition.

The isovariant theory should be considered as the equivariant theory in the unstable range. The bigger the gap is, the more equivariant the isovariant theory become. Once we reach the threshold of the big gap,  the isovariant theory and the equivariant theory become the same \cite{b}. Such viewpoint is very similar to the fact the higher connectivity makes space more homotopically stable. Indeed it is no coincidence that the dimension requirement in the Freudental suspension theorem is the same as the big gap hypothesis.

Therefore we expect to see some equivariant phenomena in isovariant theory even before reaching the big gap. For a $G$-manifold $M$, denote by $M_s=\{x\in M\colon G_x\ne \{e\}\}$ the non-free part. By considering isovariant homotopy equivalences that are already homeomorphic on the non-free part, the corresponding $G$-structure ${\mc S}_G(M,\rel M_s)$ is functorial for all equivariant maps under mild conditions.

\begin{thm}\label{th2}
Suppose $M$ and $N$ are $G$-manifolds of the same dimension. Suppose both $M_s$ and $N_s$ have codimension $\ge 3$. Then any equivariant map $f\colon M\to N$ induces a map $f_*\colon {\mc S}_G(M,\text{\rm rel } M_s) \to {\mc S}_G(N,\text{\rm rel } N_s)$. Moreover, $f_*$ depends only on the equivariant homotopy class of $f$, and $id_*=id$, $(fg)_*=f_*g_*$.
\end{thm}

By Weinberger's stratified surgery theory \cite{w}, the structure set ${\mc S}_G(M,\rel M_s)$ relative to the singular part is the fibre of the assembly map of the generalized homology theory with the $L$-spectrum (over the orbit category of group actions) as the coefficient. Such assembly map has been defined by Davis and L\"uck in \cite{dl} and many other people. In this paper, we give a more explicit explanation about the connection between the structure set relative to the singular part and the assembly map.

\section{Functoriality in Free Part}

Let $X$ be a $G$-space and let $\omega:\pi_1(X/G)\ra\{\pm 1\}$ be an orientation class. We define {\em free surgery problems} over $X$ by considering the following $G$-surgery problems with references to $X$:
\begin{equation}\label{1}
\begin{CD}
\nu_M @>\bar{f}>> \eta \\
@VVV  @VVV \\
(M,\pa M) @>f>> (N,\pa N) @>\alpha>> X
\end{CD}
\end{equation}
where
\begin{enumerate}
\item $G$ acts freely on $M$ and $N$;
\item $(\bar{f},f)$ is an equivariant degree 1 normal map;
\item $\pa f:\pa M\ra\pa N$ is a $G$-homotopy
equivalence;
\item $\omega\alpha_{*}: \pi_1(N/G)\ra\pi_1(X/G)\ra\{\pm 1\}$ is the orientation of the surgery problem.
\end{enumerate}
Denote by $L_{G,\text{free}}(X,\omega)$ the normal cobordism classes of these free surgery problems over $X$. This is an abelian group, with disjoint union as the sum and the orientation reversing as the negative operation. Moreover, the definition can be spacified, so that the surgery obstruction group is the zeroth homotopy group of the surgery obstruction space, and the higher homotopy groups are the obstructions of higher dimensional surgery problems.

We emphasis that by the very definition, $L_{G,\text{free}}(X,\omega)$ is a covariant functor of equivariant maps. This is the key motivation behind our definition.

Let $X_s=\cup_{1\neq H\sub G}X^{H}$ be the nonfree part of $X$. The following proposition is a consequence of the fact that surgery obstructions are concentrated on the 2-skeleton.  

\begin{prop}\label{p1}
Suppose that $X-X_s\sub X$ is 2-connected. Then the natural map
\[
L((X-X_s)/G,\omega)=L_G(X-X_s,\omega)
\ra L_{G,\text{\rm free}}(X,\omega)
\]
is an equivalence for surgery problems of dimension $\geq 5$.
\end{prop}

In case $X$ is a $G$-manifold, the connectivity condition is satisfied when the nonfree part has codimension $\geq 3$.

\begin{proof}

To show the map is surjective, we begin with a free surgery obstruction represented by \eqref{1}. We thicken an equivariant 2-skeleton in $N-\pa N$ to obtain a codimension 0 $G$-submanifold $(N_2,\pa N_2)$. The complement $N_{>2}=\overline{N-N_2}$ is a codimension 0 $G$-submanifold with boundary $\pa N_{>2}=\pa N_2\coprod\pa N$. Thus $N$ decomposes as $N=N_2\cup_{\pa N_2}N_{>2}$. Because $\dim N\geq 5$, we have $\pi_1N_2=\pi_1\pa N_2=\pi_1N_{>2}=\pi_1N$ by the Van-Kampen theorem. 

The 2-connectivity condition allows us to equivariantly homotope $\alpha$ such that $\alpha(N_2)\cap X_s=\emptyset$. Moreover, we may equivariantly homotope $f$ so that it is transverse to $\pa N_2$, which gives rise to a corresponding decomposition
\[
f=f_2\cup f_{>2}: M_2\cup M_{>2}\ra N_2\cup N_{>2}.
\]
These homotopies do not change the surgery obstruction class.

By $\pi_1N_2=\pi_1N$, we may apply the $\pi-\pi$ theorem to the surgery problem $f_{>2}:(M_{>2},\pa M_2)\ra (N_{>2},\pa N_2)$ relative to the homotopy equivalence $\pa f:\pa M\simeq\pa N$ on the other part of the boundary. We obtain a normal cobordism to a homotopy equivalence $f_{>2}':(M_{>2}',\pa M_2')\ra (N_{>2},\pa N_2)$. Glueing the normal cobordism to $f\times[0,1]$, we obtain a normal cobordism from $f$ to a surgery problem of the form
\[
(M'_2,\pa M'_2)\stackrel{f'_2}{\lra}
(N_2,\pa N_2)\stackrel{\alpha_2}{\lra}
X-X_s\sub X,
\]
where $\pa M'_2\simeq \pa N_2$. This proves the surjectivity.

The injectivity may be proved similarly. Given a normal cobordism in $L_{G,\text{free}}(X,\omega)$ between two surgery problems that miss $X_s$ (meaning coming from $L_G(X-X_s,\omega)$), we may take the 2-skeleton of the cobordism based on the two surgery problems. Then after a homotopy that moves the 2-skeleton away from $X_s$, the complement of a thickening of the 2-skeleton is normally cobordant to a simple homotopy equivalence by the $\pi-\pi$ theorem. We then obtain a new normal cobordism between the two origional problems that misses $X_s$. This proves the injectivity.

\end{proof}

Observe that $L_G(X-X_s,\omega)=L_G(X,\text{rel sing})$ is the surgery obstruction relative to the singular part of the group action. Moreover,  the forgetful map $L_G(X-X_s,\omega)\to L_{G,\text{free}}(X)$ is naturally defined, commuting with the assembly maps
\begin{equation}\label{2}
\begin{CD}
H(M/G,L^{-\infty}_G(\text{loc}M,\text{rel sing})) 
@>>> 
L^{-\infty}_G(M,\text{rel sing}) \\
@VV\simeq V  @VV\simeq V \\
H(M/G,L^{-\infty}_{G,\text{free}}(\text{loc}M)) @>>> L^{-\infty}_{G,\text{free}}(M)
\end{CD}
\end{equation}
By Weinberger's stratified surgery theory, the homotopy fibre of the upper arrow is the stable structure $S^{-\infty}_G(M,\text{rel sing})$. The equivariant functoriality of $L^{-\infty}_{G,\text{free}}$ implies the equivariant functoriality of the fibre of the homotopy fibre of the lower arrow. Thus we obtain the stable version of Theorem \ref{th2}.

Results similar to Theorem \ref{th2} have been obtained in special cases before. Given codimension $\geq 3$ gap condition, Cappell and Weinberger identified the normal invariants for the structure set $S_G^{-\infty}(M,\text{\rm rel sing})$ with $H(M/G,L^{-\infty}(G_x))$, simply because the free part of $\text{loc}M$ is simply connected with free isotropy group action. Then the homology $H(M/G,L^{-\infty}(G_x))$ appears to be equivariant, simply because the inclusion $G_x\sub G_{f(x)}$ induces a map on the surgery obstructions. However, there are two technical problems with this argument. One is the choice among many conjugate isotropy subgroups for a point in $M/G$. A more serious problem is compatibility between the induced maps $L(G_x)\to L(G_{f(x)})$ and the map between the homologies. 

Now we discuss the issue of whitehead torsions. The idea for deriving the functorial properties of rel sing equivariant surgery obstructions by introducing $L_{G,\text{free}}$ can also be applied to the whitehead torsions. Define $Wh^{PL}_{G,\text{free}}(X)$ by considering finite relatively free $G$-complexes $(Y,X)$ (i.e., $G$ acts freely on $Y-X$) that deformation retracts to $X$, modulo the equivalence relation given by expansion and collapsing through such complexes. The proof of Proposition \ref{p1} is based on the fact that surgery obstructions are concentrated on the 2-skeleton, a property shared by the whitehead torsions. Thus it is not surprising we should have the following result.

\begin{prop}\label{p3} 
Suppose that $X-X_s\sub X$ is 2-connected. Then the natural map
\begin{equation}\label{4}
Wh^{PL}((X-X_s)/G)=Wh^{PL}_G(X-X_s)
\ra Wh^{PL}_{G,\text{\rm free}}(X)
\end{equation} 
is an equivalence.
\end{prop}

The proposition may also be proved for the other torsion (finiteness and negative $K$-theory) invariants of lower dimensions. This is the consequence of Proposition \ref{p3} and the fact that the other torsion invariants are whitehead torsions over $X\times T^{i}$ that are invariant under all the transfers induced by various self coverings of $T^i$.

Proposition \ref{p3} (and the similar version for the other torsion) has two implications. First, Proposition \ref{p1} will also hold for the stable surgery obstructions $L^{-\infty}$, using Rothenberg sequence, so that the stable version of Theorem \ref{th2} is indeed proved. Second, Theorem \ref{th2} can be similarly proved for $Wh^{Top,\leq 0}$. This can then be used to destablize the stable version of Theorem \ref{th2}. 

\begin{proof}
Given any $k>1$, it was proved in the classical simple homotopy theory that a deformation retract over $X$ is equivalent through expansions and collapsings to a deformation retract with only $k$- and $(k+1)$-cells (subsequently called ($k$-$(k+1)$)-complex). With the same argument, an element in $Wh^{PL}_{G,\text{free}}(X)$ may be represented by a relatively free (2-3)-complex $G$-deformation retract $(Y,X)$. Since the (2-3)-complex is attached to $X$ along a 2-skeleton of $X$, the 2-connectivity assumption implies that the attaching map may be homotoped to miss $X_s$. The homotopy does not change the element in the whitehead torsion group. Thus we proved that the element comes from $Wh^{PL}_G(X-X_s)$.

For the injectivity, we note that the 2-connectivity implies that $\pi_1(X)=\pi_1(X-X_s)$. Therefore we have corresponding inclusion of universal covers $\widetilde{X-X_s}\sub\tilde{X}$. The group $\tilde{G}=\pi_1((X-X_s)/G)$ acts on $\tilde{X}$ and the action is free when restricted to $\widetilde{X-X_s}$. Now given a relatively free (2-3)-complex $G$-deformation retract $(Y,X)$, we take the universal cover to obtain a relatively free (2-3)-complex $\tilde{G}$-deformation retract $(\tilde{Y},\tilde{X})$. The boundary map $C_{3}(\tilde{Y},\tilde{X})\ra C_2(\tilde{Y},\tilde{X})$ is an isomorphsm of free ${\bf Z}\tilde{G}$-modules, which we take to define an algebraic torsion invariant $\tau(Y,X)\in Wh({\bf Z}\tilde{G})$. As in the classical case, the map is independent of the choice of the representative $(Y,X)$, so that $\tau$ is a map $Wh^{PL}_{G,\text{free}}(X)\ra Wh({\bf Z}\tilde{G})$.

What we did is simply to extend the definition of the algebraic whitehead torsion from $Wh^{PL}((X-X_s)/G)$ to $Wh^{PL}_{G,\text{free}}(X)$. Consequently, the composition of the natural map \eqref{4} with $\tau$ gives rise to the classical isomorphism $Wh^{PL}((X-X_s)/G)\cong Wh(({\bf Z}\tilde{G})$. In particular, we proved that the natural map \eqref{4} is injective.

\end{proof}

\section{The Assembly Map and the Isomorphism Conjecture}

We explain that the assembly map
\[
H(M/G;L^{-\infty}_G(\text{loc}M,\text{rel sing})) 
\ra
L^{-\infty}_G(M,\text{rel sing})
\]
that computes the stable structure set $S^{-\infty}_G(M,\text{rel sing})$ is the same assembly map for the generalized homology theory with the $L$-spectrum (over the orbit category of group actions) as the coefficient, as defined by Davis and L\"uck in \cite{dl}, for example. In fact, Davis and L\"uck showed that the assembly map has axiomatic characterization. Since Weinberger only used such axiomatic properties in establishing the stratified surgery theory in \cite{w}, the two assemblies are really the same. 

Since the argument in \cite{w} is somewhat implicit, we will give here a more explicit explanation in a simplest case (but contains all the key ingredients). 

Let us consider the case that a semi-free action of $G$ on a manifold $M$. We assume that near the fixed set $M^G$, the action is locally linear and is a block bundle over $M^G$. We denote the block bundle by $E_M\to M^G$, which has the unit ball $D(\rho)$ of a $G$-representation $\rho$ as fibre. Then consider the structure set $S_G(M,\text{rel sing})=S_G(M,\text{rel }M^G)$, consisting of equivalent classes of maps $f\colon N\to M$, where 
\begin{enumerate}
\item $N$ is a semi-free $G$-manifold, such that the action is locally linear and is a block bundle over $N^G$.
\item $f$ is an isovariant simple homotopy equivalence, such that the restriction $f^G\colon N^G\to M^G$ is a homeomorphism. 
\item $f$ is a homotopically transverse isovariant homotopy equivalence: $f$ restricts to a blockwise simple homotopy equivalence $f_E\colon E_N\to E_M$ on the block bundle neighborhoods of the fixed points. 
\end{enumerate}
The equivalence between various $f$ means homotopy equivalence. It is also possible to consider the usual (not necessarily simple) homotopy equivalences, but the block bundle under consideration needs to be relaxed to ``$h$-block bundle'', which means that a block over a simplex $\Delta$ is only homotopy equivalent to (instead of homeomprhic to) $D(\rho)\times \Delta$.

Let $M'=\overline{M-E_M}$ be the complement of the block bundle neighborhood. We have the spherical bundle $\pa E_M\to M^G$ (with the unit sphere $S(\rho)$ of $\rho$ as fibre) associated to the block bundle $E_M\to M^G$. Then $M'$ is a free $G$-manifold with $\pa M'=\pa E_M$. The structure set $S_G(M,\text{rel }M^G)$ can be ``spacified'' and fits into a fibration
\[
S_G(M',\text{rel }\pa M')\to S_G(M,\text{rel }M^G)\to S_G(E_M\to M^G,\text{rel }M^G).
\]
By taking the quotient of the action, the structure set on the left becomes the usual structure set
\[
S_G(M',\text{rel }\pa M')
=S(M'/G,\text{rel }\pa M'/G).
\]
By coning, the structure set on the right is the same as the structure set for the associated sphere block bundle
\[
S_G(E_M\to M^G,\text{rel }M^G)=S_G(\pa E_M\to M^G,\text{rel }M^G).
\]
This is the equivalence classes of blockwise simple equivariant homotopy equivalences
\begin{equation}\label{3}
\begin{CD}
W @>>> \pa E_M \\
@VVV  @VVV \\
M^G @= M^G
\end{CD}
\end{equation}
The quotient by $G$ are blockwise simple homotopy equivalences
\begin{equation}\label{6}
\begin{CD}
V @>>> \pa E_M/G \\
@VVV  @VVV \\
M^G @= M^G
\end{CD}
\end{equation}
where $\pa E_M/G\to M^G$ is a block bundle with manifold (a spherical space form for the semi-free action) $F=S(\rho)/G$ as fibre. Since \eqref{3} can be recovered from \eqref{6} by pullback (or taking cover in case $G$ is finite), we see that
\[
S_G(\pa E_M\to M^G,\text{rel }M^G)
=S(\pa E_M/G\to M^G,\text{rel }M^G),
\]
and the structure set $S_G(M,\text{rel }M^G)$ fits into a fibration
\begin{equation}\label{7}
S(M'/G,\text{rel }\pa M'/G)\to S_G(M,\text{rel }M^G)\to S(\pa E_M/G\to M^G,\text{rel }M^G),
\end{equation}
in which the other two terms can be computed by non-equivariant surgery theory. We will show that the fibration is consistent with the assembly map for the generalized homology theory with the $L$-spectrum.

The assembly map for the generalized homology theory fits into a commutative diagram
\begin{equation}\label{5}
\begin{CD}
H(M'/G; L_G(\text{loc}M'))  
@>>>  L_G(M')  \\
@VVV  @VVV \\
H(M/G;L_G(\text{loc}M,\text{rel sing})) 
@>>> L_G(M,\text{rel sing}) \\
@VVV   \\
H(M/G,M'/G;L_G(\text{loc}M,\text{rel sing})) 
\end{CD}
\end{equation}
The two horizontal maps are the assembly maps for the generalized homology theories. The left is the fibration for the generalized homology theory of the pair $(M/G,M'/G)$. Moreover, the coefficient $L_G(\text{loc}M,\text{rel sing})=L_G(\text{loc}M')$ on the free part $M'$ of the action. 

Since $M'$ is free part of the group action, at any point in $M'$, we have $\text{loc}M=G\times D$ for a small disk $D$, and
\[
L_G(\text{loc}M,\text{rel sing})
=L_G(G\times D)
=L(D)
=L(*)
\]
is the basic $L$-spectrum. Therefore the assembly map at the top of \eqref{5} is the classical assembly map for the classical surgery theory. The fibre of the assembly map is simply the structure set $S(M'/G,\text{rel }\pa M'/G)$.

The surgery theory we tried to establish says that the fibre of the assembly map at the middle of \eqref{5} is $S_G(M,\text{rel }M^G)$. Since 
\[
L_G(M,\text{rel sing})
=L_G(M,\text{rel }M^G)
=L_G(M'),
\]
the problem is reduced to showing that
\[
S(M'/G,\text{rel }\pa M'/G)
\to S_G(M,\text{rel }M^G)
\to H(M/G,M'/G;L_G(\text{loc}M,\text{rel sing})) 
\]
is a fibration. Compared with the fibration \eqref{7}, we simply need to establish
\[
H(M/G,M'/G;L_G(\text{loc}M,\text{rel sing})) 
=S(\pa E_M/G\to M^G,\text{rel }M^G).
\]
By excision for the generalized homology theory, this is the same as
\[
H(E_M/G,\pa E_M/G;L_G(\text{loc}M,\text{rel sing})) 
=S(\pa E_M/G\to M^G,\text{rel }M^G).
\]
Then by the fibration for the generalized homology theory of the pair $(E_M/G,\pa E_M/G)$, this means that
\[
H(\pa E_M/G;L_G(\text{loc}M,\text{rel sing}))
\to H(E_M/G;L_G(\text{loc}M,\text{rel sing}))
\to S(\pa E_M/G\to M^G,\text{rel }M^G)
\]
is a fibration. By the homotopy invariance of the generalized homology theory, the problem is finally reduced to the establishment of the fibration
\begin{equation}\label{8}
H(\pa E_M/G;L_G(\text{loc}M,\text{rel sing}))
\to H(M^G;L_G(\text{loc}M,\text{rel sing})) 
\to S(\pa E_M/G\to M^G,\text{rel }M^G).
\end{equation}

To avoid complication in notation, we further assume the block bundle structure in the neighborhood of the fixed point to be trivial. In other words, we assume the block bundle $\pa E_M/G\to M^G$ is simply the projection $F\times M^G\to M^G$, with $F=S(\rho)/G$. Under such assumption, the blockwise simple homotopy equivalence over $M^G$ simply means a consistent collection of simple homotopy equivalences indexed by the simplices in $M^G$. Therefore in terms of the simplicial space, we have
\[
S(\pa E_M/G\to M^G,\text{rel }M^G)
=S(F\times M^G\to M^G,\text{rel }M^G)
=Maps[M^G,S(F)]
\]
The classical surgery fibration
\[
S(F)\to Maps[F,L(*)]\to L(F)
\]
then induces a fibration
\[
Maps[M^G,S(F)]
\to Maps[M^G,Maps[F,L(*)]]
\to Maps[M^G,L(F)].
\]
Using $Maps[F\times M^G,L(*)]]=Maps[M^G,Maps[F,L(*)]]$, this is the same as
\[
S(\pa E_M/G\to M^G,\text{rel }M^G)
\to Maps[\pa E_M/G, L(*)]
\to Maps[M^G,L(F)].
\]
The space of maps from a space to an $L$-spectrum is simply the generalized cohomology theory of the space with the $L$-spectrum as the coefficient. In our case, the spaces $\pa E_M/G$ and $M^G$ are manifolds, and generalized Poincar\'e duality holds. The result is the following fibration
\begin{equation}\label{9}
S(\pa E_M/G\to M^G,\text{rel }M^G)
\to H(\pa E_M/G;L(*))
\to H(M^G;L(F)).
\end{equation}
Since $G$ acts freely on $\pa E_M$, we have
\begin{equation}\label{a}
L_G(\text{loc}M,\text{rel sing})
=L(*)
\end{equation}
along $\pa E_M$. On the other hand, at any point on $M^G$, we have $\text{loc} M=cF\times D$, where $cF$ is the cone on $F=S(\rho)/G$, and $D$ is a disk in $M^G$. Therefore we have
\begin{equation}\label{b}
L_G(\text{loc}M,\text{rel sing})
=L(cF\times D,\text{rel }D)
=L(F)
\end{equation}
along $M^G$. Therefore the fibration \eqref{9} is
\begin{equation}\label{10}
S(\pa E_M/G\to M^G,\text{rel }M^G)
\to H(\pa E_M/G;L_G(\text{loc}M,\text{rel sing}))
\to H(M^G;L_G(\text{loc}M,\text{rel sing})).
\end{equation}
This apparent difference from \eqref{8} is only due to the fact that we did not keep track of the dimensions of the spaces. For example, $\text{loc} M$ is really one dimension higher than the dimension of $\pa E_M/G$. If we really keep track of dimensions, then \eqref{a} and \eqref{b} should really be
\[
L_G(\text{loc}M,\text{rel sing})
=\Omega L(*),\s
L_G(\text{loc}M,\text{rel sing})
=\Omega L(F),
\]
so that the fibration \eqref{10} is really
\[
S(\pa E_M/G\to M^G,\text{rel }M^G)
\to \Sigma H(\pa E_M/G;L_G(\text{loc}M,\text{rel sing}))
\to \Sigma H(M^G;L_G(\text{loc}M,\text{rel sing})).
\]
This is the same fibration as \eqref{8}.

There is only technical problem in the explanation above. The block bundle structure exists in the PL category, while the basic building block of the theory is the classical surgery fibration 
\[
S(M)\to Maps[M,L(*)]\to L(M),
\]
which only applies to the topological manifolds $M$. The fibration needs to be modified a little bit in the low dimension in order to be applied to PL manifolds. But the modification destroys the fibration structure. 

What solves the technical problem is the stablization process. Specifically, the analysis of the existence of PL structure on locally triangulable spaces by Anderson and Hsiang \cite{ah1,ah2} and Quinn \cite{q1} implies that, for a manifold homotopically stratified space $X$, there is a big $k$ depending only on the dimensions of the strata of $X$, such that $X\times {\bb R}^k$ is a manifolds ${\bb R}^k$-controlled geometrically stratified space. Here geometrically stratified means block bundle like structure for the neighborhood of strata, and homotopically stratified (see \cite{q2}) means that the homotopy properties of geometrical stratification are satisfied. For sufficiently big $k$, we then have the stable structure 
\[
S^{-\infty}_G(X)=S_G(X\times {\bb R}^k,\text{ controlled over }{\bb R}^k),
\]
and the stable surgery obstruction 
\[
L^{-\infty}_G(X)=L_G(X\times {\bb R}^k,\text{ controlled over }{\bb R}^k).
\]
Moreover, we have the stablization maps
\[
\times {\bb R}^k\colon 
S_G(X)\to S^{-\infty}_G(X),\s
L_G(X)\to L^{-\infty}_G(X).
\]
Since we have block bundle like structure after stablization, the argument above on the stratified surgery fibration is valid in the stable range.

\section{Comments}

At the first sight, the functoriality is rather mysterious, even for the non-equivariant homotopy classification. It is not at all obvious that given any map $f\colon M\to N$, one can naturally associate to a homotopy equivalence $M'\simeq M$ another homotopy equivalence $N'\simeq N$. The following is the explanation in the non-equivariant case, given on page 82 of \cite{w}: Embed $N$ into a sphere $S^n$ of sufficiently big dimension $n$. The regular neighborhood ${\mc N}_n$ of $N$ in the sphere is an $n$-dimensional manifold with boundary $\pa {\mc N}_n$. Denote by $S_n(N)$ the structure set $S({\mc N}_n,\rel \pa {\mc N}_n)$ of ${\mc N}_n$ that are already homeomorphic on the boundary. Now because the dimension $n$ is big, it is possible to homotopically perturb the map $f\colon M\to N$ into an embedding $M\to S^n$, such that a regular neighborhood ${\mc M}_n$ of $M$ in $S^n$ can be found to lie completely inside the interior of ${\mc N}_n$. Since the structures in $S_n(M)$ are already homeomorphic on the boundary $\pa {\mc M}_n$, we have a map 
\[
\cup\overline{{\mc M}_n-{\mc N}_n}\colon
S_n(M)=S({\mc M}_n,\rel \pa {\mc M}_n)\to S_n(N)=S({\mc N}_n,\rel \pa {\mc N}_n)
\]
obtained by glueing the complement of ${\mc M}_n$ on ${\mc N}_n$ along $\pa {\mc M}_n$. Now suppose $M$, $N$ are orientable topological manifolds, and $n$ differs from the dimension of $M$ and $N$ by multiples of $4$, then by the periodicity \cite{cw, ks}, we get
\[
f_*\colon S(M)=S_n(M)\to S(N)=S_n(N).
\]

We note that $S_n(X)$ can be introduced for any finite complex $X$ and sufficiently big $n$. The periodicity tells us that $S_n(X)$ depends only on $n$ mod $4$ and may be considered as the ``stable structure set'' of $X$. The fact that $S_n(X)$ can be defined in the very general context is the key reason that the functoriality exists in the stable range. The periodicity then allows us to ``destablize'' to the actual structure of $X$ in case $X$ is a manifold. The result is the ``unstable'' functoriality.

Now let us examine the idea in the equivariant case. The sphere $S^n$ may be replaced by a sufficiently big representation (say, multiples of the regular representation) $V$. Then any fairly nice $G$-space $X$ may be embedded into sufficiently big representation. Equivariant regular neighborhoods may be chosen, the stable structure set $S_V(X)$ may be introduced, and we have the stable equivariant functoriality. Using the equivariant periodicity \cite{wy1, wy2, y1, y2}, we get $S_{V\oplus W\oplus W}(X)=S_V(X)$ for sufficiently sophisticated $V$ and any complex $G$-representation $W$. However, in case $X=M$ is a $G$-manifold, we may run into problem when we try to ``destablize''. The problem is that when we try to apply the equivariant periodicity to obtain $S_G(M)=S_G({\mc M}_V,\rel \pa {\mc M}_V)=S_V(M)$, we need to assume that there is a one-to-one correspondence between the fixed point components of $M$ and ${\mc M}_V$ (and the additional problem that the normal representation of $M$ in $V$ is a ``periodicity representation''). This difficulty prevents us to use the embedding idea to get the functoriality in general. 

So we will appeal to some machinery that is less geometrical in order to get the functoriality. In the classical case, this means that we consider the surgery fibration
\[
S(M)\to Maps[M,G/Top]\to L(M),
\]
whose associated long exact sequence of homotopy groups is the usual surgery exact sequence \cite{q}. Then Siebemann's periodicity \cite{ks} implies that $L(e)={\bb Z}\times G/Top$, and the space is actually an infinite loop space. Therefore up to an additional copy of ${\bb Z}$, $Maps[M,G/Top]$ is the generalized cohomology theory with the spectrum $L(e)$ as the coefficient. Moreover, the generalized Poincar\'e duality can be applied to the cohomology theory and turns $Maps[M,L(e)]$ into a homology theory. The map from the homology theory to the surgery obstruction $L(M)$ is then the classical assembly map. Since both the homology theory and the surgery obstruction are convariantly functorial, we conclude that the structure set $S(M)$ is also convariantly functorial. 

In the equivariant case, Weinberger \cite{w} showed that, after stablization, the isovariant structure set is still the fibre of the assembly map of a sophisticated homology theory with surgery obstruction of the local isovariant structures as the coefficients. This allows us to extend the classical functoriality to the isovariant case.

\end{document}